\documentclass{article}
\usepackage{amsmath,amssymb,pxfonts}
\newtheorem{lem}{Lemma}[section]
\newtheorem{thm}{Theorem}[section]
\newtheorem{cor}{Corollary}[section]
\newtheorem{rem}{Remark}[section]
\title{EM algorithm for  stochastic hybrid systems}
\author{Masaaki Fukasawa
\vspace*{0,5cm}\\
{\small Graduate School of Engineering Science}
\\ 
{\small Osaka University}
\vspace*{0,5cm}\\
{\small 560-8531, Osaka, Japan}\\
{\small fukasawa@sigmath.es.osaka-u.ac.jp
}}
\date{}
\begin{document}
\maketitle
\begin{abstract}
A stochastic hybrid system, also known as a switching diffusion, is a
 continuous-time Markov process with state space consisting of 
 discrete and continuous parts.
We consider parametric estimation of the Q matrix for the discrete state
 transitions and of the drift coefficient for the diffusion part.
First, we derive the likelihood function under the complete observation
 of a sample path in continuous-time. Then,
extending a finite-dimensional filter for hidden Markov models developed by Elliott 
et al.~({\it Hidden Markov Models}, Springer, 1995) to stochastic hybrid systems, we derive the likelihood function and the EM algorithm
under a partial observation where
the continuous state is monitored continuously in time, while the discrete
 state is unobserved.
\\


\end{abstract}
\section{Introduction}
A stochastic hybrid system (SHS, hereafter),
also known as a switching diffusion~\cite{GAM},
is a continuous-time Markov process $Z$ with
state space $S = \{e_1,\dots,e_k\}\times \mathbb{R}^d$
consisting of both discrete and continuous parts, namely,
$\{e_1,\dots,e_k\}$ and $\mathbb{R}^d$ respectively.
The elements $\{e_i\}$ are, without loss of generality, specified as
the standard basis of $\mathbb{R}^k$ in this article.
Denoting by $\langle \cdot, \cdot \rangle$ the inner product of
$\mathbb{R}^k$ or $\mathbb{R}^d$,
$\langle e_i,e_j \rangle = \delta_{ij}$, where $\delta_{ij}$ is 
Kronecker's delta.
The discrete part of $Z$, denoted by $X$, is a step process  with 
state space $\{e_1,\dots,e_k\}$ and 
 ``$Q$ matrix''
of the form $Q(Y_t) = [q_{ij}(Y_t)]$, where $Y$ is the continuous part of
$Z$. In other words, 
\begin{equation}\label{tr.prob}
 P(X_{t+h} = e_j | X_t = e_i, Y_t = y) = (\delta_{ji} + q_{ji}(y))h + o(h) 
\end{equation}
as $h \to 0$.
Here, $Q(y) = [q_{ij}(y)]$ is a $Q$ matrix for each $y$, that is,
$q_{ji}(y) \geq 0$ for $j \neq i$ and $\sum_{j \neq i} q_{ji}(y) = -q_{ii}(y)$ for 
all $y \in \mathbb{R}^d$.
The continuous part $Y$ is 
defined as the solution of a stochastic differential equation
\begin{equation*}
\dot{Y} = \mu(X,Y) + \epsilon \dot{W}
\end{equation*}
 on $\mathbb{R}^d$ for some $\mathbb{R}^d$-valued function $\mu$,
where $\epsilon > 0$ and 
$\dot{W}$ is a $d$ dimensional Gaussian white noise.
The generator $\mathcal{L}$ of the Markov process $Z = (X,Y)$ is given by
\begin{equation}\label{gen}
\mathcal{L}f(e_i,y) =  \langle \mu,  \nabla_yf \rangle (e_i,y)+ \frac{1}{2}\epsilon^2
\Delta_y f(e_i,y) + \sum_{j=1}^k (f(e_j,y) - f(e_i,y))\langle e_j, Q(y)e_i\rangle,
\end{equation}
where $\nabla_y$ and $\Delta_y$ refer to the gradient and Laplacian
operators respectively with respect to the variable $y$.

There is a huge amount of literature on the analysis and applications of  
SHS. See e.g., \cite{GAM, H, TSS, Rev, ME, YZ} and the references therein.
The author's motivation to study the SHS is its potential application to 
the analysis of certain single-molecule dynamics which have several hidden
states with the switching rates depending on the (observable)
position of the molecule~\cite{Iwaki}. 
In this article, our aim is to construct an estimation algorithm
for the Q matrix 
$Q(y) = Q^{\theta}(y)$
and of the drift coefficient $\mu(z) = \mu^{\theta}(z)$
based on a continuous  path of 
$Y$ with $X$  unobserved.

Our problem is closely related to filtering. Indeed,
when $Q$  does not depend on $y$, the system is a hidden
Markov model,
for which finite dimensional filters are given by 
 Wonham~\cite{Wonham} and Elliott et
al.~\cite{El}.
Dembo and Zeitouni~\cite{DZ} combined the filtering and the EM algorithm
for parametric estimation under partial observations.
Although the general case of our interest is beyond the framework 
of \cite{DZ, El}, we observe in this article that the arguments in \cite{DZ,El}
can be extended.

A general theory of filtering, estimation and related stochastic analysis
is available in Liptser and Shiryaev~\cite{LS1,LS2}.
In general, an optimal filtering equation is 
infinite dimensional, which hampers
its direct use in practice. 
Finite dimensional filters have been obtained in some cases (see Bain
and Crisan~\cite{BC}) including the hidden Markov model as mentioned
above,
of which the essential property is that the state space of the hidden
process is finite.
The literature of filtering is huge. Nevertheless, a filtering problem
considers a signal (hidden) process 
of which the law is determined irrespectively of an observable process.
On the other hand,
our interest is  in the case where a hidden process and an observable
process are coupled.

In Section~2, we describe the basic properties of SHS as the solution of a
martingale problem. In Section~3, we derive the likelihood function
under complete observations of both $X$ and $Y$ on a time interval $[0,T]$.
In Section~4, we consider the case where the discrete part $X$ is
unobservable, and construct a finite dimensional filter extending \cite{El}.
In Section~5, using the finite dimensional filter, 
we derive the likelihood function under the partial
observation.
In Section~6, extending \cite{DZ,El},
we construct the EM algorithm for 
an exponential family
 under the partial  observation.

\section{A construction as a weak solution}
Here we construct a SHS as a weak solution, that is, we construct a
distribution on the path space
$D([0,T];S)$ which is a solution of the martingale problem with the
generator (\ref{gen}). 
\\

A direct application of Theorem (5.2) of Stroock~\cite{Str} 
provides the following.
\begin{thm}\label{thm0}
 Let $\mu$ be a bounded Borel function and 
$q_{ij}$, $1 \leq i,j \leq k$ be bounded continuous
 functions. Then, for any $z \in S$,
there exists a unique probability measure $P_z$ on
$D([0,T];S)$ such that $Z_0 = z$ and 
\begin{equation*}
 f(Z_t) - \int_0^t \mathcal{L}f(Z_s)\mathrm{d}s
\end{equation*}
is a martingale under $P_z$ for any $f \in
 C^{0,\infty}_0(\{e_1,\dots,e_k\}\times \mathbb{R}^d)$,
where $Z: t \mapsto Z_t$ is the canonical map on $D([0,T];S)$.
Moreover, $Z$ is a strong Markov process with $\{P_z\}_{z \in S}$.
\end{thm}
The uniqueness part of Theorem~\ref{thm0} is important in this
article.
For the existence, we give below an explicit construction, which plays
a key role 
to solve 
a filtering problem later.\\

First, we construct a SHS with $\mu = 0$ in a pathwise manner.
Without loss of generality, assume $\epsilon = 1$.
Note that $Y$ is then a $d$ dimensional Brownian motion.
Let $(\Omega,\mathcal{F},P^0)$ be a probability space on which a $d$
dimensional Brownian motion $Y$ and an i.i.d.~sequence of exponential random
variables $\{E_n\}$ that is independent of $Y$ are defined.
Conditionally on $Y = \{Y_t\}_{t\geq 0}$, a time-inhomogeneous
continuous-time Markov chain $X$ with (\ref{tr.prob})
 is defined using the exponential variables. More specifically,
given $X_0 = e_i$, let
\begin{equation*}
\tau_1 = \min_{1 \leq j \leq k} \tau_1^j, \ \ 
 \tau_1^j = \inf\left\{t > 0 ; \int_0^t q_{ji}(Y_s)\mathrm{d}s > E_j\right\}
\end{equation*}
and $X_t = X_0$ for $0 \leq t < \tau_1$, $X_{\tau_1} = e_J$ with 
$J = \mathrm{argmin} \tau_1^j$. The construction goes in a recursive
manner;
given $X_{\tau_n} = e_i$,
let
\begin{equation*}
\tau_{n+1} = \min_{1 \leq j \leq k} \tau_{n+1}^j, \ \ 
 \tau_{n+1}^j = \inf\left\{t > \tau_n ; \int_{\tau_n}^t q_{ji}(Y_s)\mathrm{d}s >
		     E_{nk + j}\right\}
\end{equation*}
and $X_t = X_{\tau_n}$ for $\tau_n \leq t < \tau_{n+1}$, 
$X_{\tau_{n+1}} = e_J$ with 
$J = \mathrm{argmin} \tau_{n+1}^j$. 
Properties of the exponential distribution verifies the following lemma.
\begin{lem}\label{lem1}
Assume $q_{ij}(y)$ is bounded and continuous in $y \in \mathbb{R}^d$
for each $(i,j)$. Then,
\begin{equation}\label{tr.prob2}
 P^0(X_{t+h} = e_j | X_t = e_i, Y) = (\delta_{ji} + 
q_{ji}(Y_t))h + o(h)
\end{equation}
and (\ref{tr.prob}) with $P = P^0$.
\end{lem}

By It\^o's formula, for any $f \in C^{0,2}_b(\{e_1,\dots,e_k\}\times
\mathbb{R}^d)$, we have
\begin{equation*}
\begin{split}
 f(X_{t+h},Y_{t+h}) = f(X_t,Y_t) + &
  \int_t^{t+h}\langle \nabla_yf(X_s,Y_s), \mathrm{d}Y_s\rangle
+ \frac{1}{2}\epsilon^2 \int_t^{t+h}\Delta_y f(X_s,Y_s)\mathrm{d}s \\&
+ \sum_{t < s \leq t+h} (f(X_s,Y_s)-f(X_{s-},Y_s)),
\end{split}
\end{equation*}
from which together with (\ref{tr.prob}) it follows
\begin{equation*}
 \lim_{h \to 0} \frac{E^0[f(X_{t+h},Y_{t+h})|X_t = e_i,Y_t=y] -
  f(e_i,y)}{h}
= \mathcal{L}^0f(e_i,y),
\end{equation*}
where $E^0$ is the expectation under $P^0$ and
$\mathcal{L}^0f = \mathcal{L}f$ with $\mu = 0$ in (\ref{gen}).
 Note only this, we have also that 
\begin{equation*}
U^{0,f}_t :=  f(X_t,Y_t) - \int_0^t \mathcal{L}^0f(X_s,Y_s)\mathrm{d}s
\end{equation*}
is a martingale with respect to the filtration  $\{\mathcal{F}_t\}$ generated by $Z = (X,Y)$.
Even more importantly, Lemma~\ref{lem1} implies the following.
\begin{lem}\label{lem2}
Under the same conditions of Lemma~\ref{lem1}, for any 
function $g$ on $\{e_1,\dots,e_k\}$,
\begin{equation*}
V^{0,g}_t :=  g(X_t) - \int_0^t \mathcal{L}^0g(X_s)\mathrm{d}s
\end{equation*}
is a martingale with respect to the natural filtration of $X$
under the conditional probability measure given $Y$,
where $\mathcal{L}^0g = \mathcal{L}^0f$ with $f(x,y) = g(x)$.
In particular,
\begin{equation*}
 X_t - \int_0^t Q(Y_s)X_s\mathrm{d}s
\end{equation*}
is a martingale under the conditional probability measure $P^0(\cdot | Y)$.
\end{lem}
\begin{rem}
 \upshape
Lemma~\ref{lem2} will be the key to extend a finite dimensional filter
developed in Chapter~8 of Elliott et al.~\cite{El}.
The measure $P^0$ corresponds to $\bar{P}$ in 
Chapter~8 of \cite{El}, under which $X$ and $Y$ are independent.
In our framework, when $Q$ depends on $y$, they are not independent
anymore but still,
$Y$ is a Brownian motion.
 This, together with the martingale property in
Lemma~\ref{lem2}, enables us to compute conditional expectations in
 Theorem~\ref{thm3} below.
\end{rem}

Now we construct a SHS for a general bounded Borel function $\mu$.
Let
\begin{equation}\label{Lambda}
 \Lambda_T = \exp\left\{
\frac{1}{\epsilon^2}\int_0^T \langle \mu(X_t,Y_t), \mathrm{d}Y_t \rangle 
- \frac{1}{2\epsilon^2}\int_0^T |\mu(X_t,Y_t)|^2
\mathrm{d}t
\right\}.
\end{equation}
By the boundedness of $\mu$, Novikov's conditions is satisfied and so,
$\Lambda$ is an $\{\mathcal{F}_t\}$-martingale under $P^0$. 
Therefore,
\begin{equation*}
 \frac{\mathrm{d}P}{\mathrm{d}P^0} = \Lambda_T
\end{equation*}
defines a probability space $(\Omega,\mathcal{F}_T,P)$.
\begin{thm}\label{thm1}
Let $Q = [q_{ij}]$ be a Q matrix-valued bounded continuous function
and $\mu$ be an $\mathbb{R}^d$-valued bounded Borel function.
 Under $P$,
$Z = (X,Y)$ is a Markov process with generator (\ref{gen}).
Further for any $f \in C^{0,2}_b(\{e_1,\dots,e_k\}\times
\mathbb{R}^d)$,  
\begin{equation*}
U^f_t :=  f(Z_t) - \int_0^t \mathcal{L}f(Z_s)\mathrm{d}s
\end{equation*}
is an $\{\mathcal{F}_t\}$ martingale.
\end{thm}
{\it Proof:} By the Bayes formula,
\begin{equation*}
 E[f(Z_{t+h})|\mathcal{F}_t] = \frac{E^0[\Lambda_T
  f(Z_{t+h})|\mathcal{F}_t]}{E^0[\Lambda_T|\mathcal{F}_t]}
= E^0\left[
\frac{\Lambda_{t+h}}{\Lambda_t} f(Z_{t+h})|\mathcal{F}_t
\right].
\end{equation*}
Since
\begin{equation*}
 \frac{\Lambda_{t+h}}{\Lambda_t} = \exp\left\{
\frac{1}{\epsilon^2}\int_{t}^{t+h} \langle \mu(Z_s), \mathrm{d}Y_s \rangle
- \frac{1}{2\epsilon^2}\int_t^{t+h} |\mu(Z_s)|^2
\mathrm{d}s
\right\}
\end{equation*}
and $Z$ is Markov under $P^0$, 
$E[f(Z_{t+h})|\mathcal{F}_t] = E[f(Z_{t+h})|Z_t]$, meaning that 
it is Markov under $P$ as well.
By It\^o's formula, 
\begin{equation*}
 \mathrm{d}\Lambda_t = \frac{1}{\epsilon^2}\Lambda_t \mu(Z_t)\mathrm{d}Y_t
\end{equation*}
and
\begin{equation*}
\begin{split}
 \Lambda_{t+h}f(Z_{t+h}) = & \  \Lambda_tf(Z_t) + 
\int_t^{t+h}f(Z_s)\mathrm{d}\Lambda_s 
+ \int_t^{t+h}\Lambda_{s-}\mathrm{d}U^{0,f}_s  
\\ & + \int_t^{t+h}\Lambda_s\mathcal{L}^0f(Z_s)\mathrm{d}s
+ \int_t^{t+h} \Lambda_s\langle \mu, \nabla_yf \rangle (Z_s)\mathrm{d}s.
\end{split}
\end{equation*}
Therefore,
\begin{equation*}
 \Lambda_{t+h}U^f_{t+h} = \Lambda_tU^f_t +
\int_t^{t+h}U^f_s\mathrm{d}\Lambda_s 
+ \int_t^{t+h}\Lambda_{s-}\mathrm{d}U^{0,f}_s,
\end{equation*}
meaning that $\Lambda U^f$ is a martingale under $P^0$.
The Bayes formula then implies that $U^f$ is a martingale under $P$.
In particular, the generator is given by $\mathcal{L}$.
\hfill////

\begin{cor}\label{cor1}
Under the same condition of Theorem~\ref{thm1},
\begin{equation*}
V_t : =  X_t - \int_0^t Q(Y_s)X_s\mathrm{d}s
\end{equation*}
is an $\{\mathcal{F}_t\}$-martingale.
\end{cor}

By the uniqueness result of Theorem~\ref{thm0},
the law of $Z$ under $P$ coincides with $P_z$ with $z = Z_0$.

\begin{rem}
\upshape
An extension to the case where $\epsilon$ depends on $Y$ is
 straightforward under the nondegeneracy of $\epsilon$ and
the well-posedness of the SDE
\begin{equation*}
 \mathrm{d}Y_t = \epsilon(Y,t)\mathrm{d}W_t
\end{equation*} 
(but with more tedious expressions).
The key of our framework is the existence of an equivalent measure under
 which the law of $Y$ is uniquely determined irrespectively of $X$.
This property admits not only (\ref{tr.prob}) but also (\ref{tr.prob2}),
which as we will see is sufficient to follow the arguments in
 Elliott et al.~\cite{El} to derive a finite dimensional filter.
Therefore the case where $\epsilon$ depends on $X$ is fundamentally different.
\end{rem}

\section{The likelihood under complete observations}
Here we consider a statistical model $\{P^\theta\}_{\theta \in \Theta}$
and derive the likelihood under complete observation of a sample path
$Z = (X,Y)$ on a time
interval $[0,T]$.
For each $\theta \in \Theta$, $P^\theta$ denotes the distribution on
$D([0,T];S)$ induced by
 a Markov process $Z$ with generator
\begin{equation*}
 \mathcal{L}^\theta f(e_i,y) =  \langle \mu^\theta,  \nabla_y f \rangle (e_i,y)+ \frac{1}{2}\epsilon^2
\Delta_y f(e_i,y) + \sum_{j=1}^k (f(e_j,y) - f(e_i,y))\langle e_j,
Q^\theta(y)e_i \rangle,
\end{equation*}
where $\mu^\theta$ is a family of $\mathbb{R}^d$-valued bounded Borel
functions and $Q^\theta = [q^{\theta}_{ij}]$ is  a family of 
$Q$ matrix-valued bounded continuous functions.
Note that $\epsilon > 0$ is almost surely identified from a path of $Y$
by computing its quadratic variation. 
It is therefore assumed to be known hereafter.
With a slight abuse of notation, we will use $P^\theta$ also 
to mean $P$ (that is, the probability measure on $\Omega$, not on $D([0,T];S)$) 
when the true parameter is $\theta$ (that is, when
$P^\theta = P\circ Z^{-1}$ with $\mu = \mu^\theta$ and $Q
= Q^\theta$).
We assume that the initial distribution $P^\theta \circ Z_0^{-1}$ is 
 known and does not  depend on $\theta$.

\begin{thm}\label{thm2}
Let $\theta, \theta_0 \in \Theta$ and assume that
\begin{equation*}
y \mapsto \frac{q^{\theta}_{ij}(y)}{q^{\theta_0}_{ij}(y)}, \ \ 
y \mapsto \frac{q^{\theta_0}_{ij}(y)}{q^{\theta}_{ij}(y)} 
\end{equation*}
are bounded  for each $(i,j)$, where $0/0 = 1$.
Then, $P^\theta$ is equivalent to $P^{\theta_0}$, 
and the log likelihood
\begin{equation*}
L_T(\theta, \theta_0)
:= \log \frac{\mathrm{d}P^\theta}{\mathrm{d}P^{\theta_0}}(\{Z_t\}_{t
\in [0,T]})
 \end{equation*}
is given by
\begin{equation*}
 \begin{split}
  L_T(\theta, \theta_0) = &\ \sum_{i\neq j}\left\{ \int_0^T\log
  \frac{q^\theta_{ji}(Y_t)}{q^{\theta_0}_{ji}(Y_t)} \mathrm{d}N^{ji}_t
- \int_0^T (q^\theta_{ji}(Y_t) - q^{\theta_0}_{ji}(Y_t))\langle
  X_t,e_i\rangle \mathrm{d}t \right\}
\\& + \frac{1}{\epsilon^2}\int_0^T \langle
  \mu^\theta(Z_t)-\mu^{\theta_0}(Z_t),\mathrm{d}Y_t\rangle
- \frac{1}{2\epsilon^2} \int_0^T (|\mu^\theta(Z_t)|^2 -
  |\mu^{\theta_0}(Z_t)|^2)\mathrm{d}t,
 \end{split}
\end{equation*}
where 
 $N^{ji}$ is the counting process of the transition from $e_i$ to $e_j$:
\begin{equation}\label{nji}
 N^{ji}_t = \int_0^t \langle X_{s-}, e_i\rangle \langle
  e_j,\mathrm{d}X_s \rangle. 
\end{equation}
\end{thm}
{\it Proof: }
A general representation of the density process can be found in
Theorem III.5.19 of Jacod and Shiryaev~\cite{JS}. 
For such a specific model as ours, a direct derivation is more
elementary. It is standard but here given for the readers' convenience.
Let
\begin{equation*}
 L^{ji}_\tau = \int_0^\tau\log
  \frac{q^\theta_{ji}(Y_t)}{q^{\theta_0}_{ji}(Y_t)} \mathrm{d}N^{ji}_t
- \int_0^\tau (q^\theta_{ji}(Y_t) - q^{\theta_0}_{ji}(Y_t))\langle
  X_t,e_i\rangle \mathrm{d}t
\end{equation*}
and
\begin{equation*}
\begin{split}
 L^0_\tau & = 
\frac{1}{\epsilon^2}\int_0^\tau \langle
  \mu^\theta(Z_t)-\mu^{\theta_0}(Z_t),\mathrm{d}Y_t\rangle
- \frac{1}{2\epsilon^2} \int_0^\tau (|\mu^\theta(Z_t)|^2 -
  |\mu^{\theta_0}(Z_t)|^2)\mathrm{d}t \\
& = 
\frac{1}{\epsilon^2}\int_0^\tau \langle
  \mu^\theta(Z_t)-\mu^{\theta_0}(Z_t),\mathrm{d}Y_t - \mu^{\theta_0}(Z_t)\mathrm{d}t\rangle
- \frac{1}{2\epsilon^2} \int_0^\tau 
|\mu^\theta(Z_t) - \mu^{\theta_0}(Z_t)|^2\mathrm{d}t. 
\end{split}
\end{equation*}
By It\^o's formula,
\begin{equation*}
\begin{split}
 \exp\{L^{ji}_\tau\} &=  \ 1 - \int_0^\tau \exp\{L^{ji}_t\}
(q^\theta_{ji}(Y_t) - q^{\theta_0}_{ji}(Y_t))\langle
  X_t,e_i\rangle \mathrm{d}t + \sum_{0 < t \leq
 \tau}(\exp\{L_t^{ji}\}-\exp\{L_{t-}^{ji}\})
\\ &= 1 + \int_0^\tau \exp\{L_{t-}^{ji}\}\left(
\frac{q^{\theta}_{ji}(Y_t)}{q^{\theta_0}_{ji}(Y_t)} -1
\right) \left[
\mathrm{d}N^{ji}_t - q^{\theta_0}_{ji}(Y_t)\langle X_t,e_i \rangle \mathrm{d}t
\right]
\end{split}
\end{equation*}
and by (\ref{nji}),
\begin{equation}\label{nji2}
 \mathrm{d}N^{ji}_t - q^{\theta_0}_{ji}(Y_t)\langle X_t,e_i \rangle
  \mathrm{d}t
= \langle X_{t-},e_i\rangle \langle e_j, \mathrm{d}X_t - Q^{\theta_0}(Y_t)X_t
\mathrm{d}t \rangle.
\end{equation}
Therefore, by Corollary~\ref{cor1}, 
$\exp\{L^{ji}\}$ and $\exp\{L^0\}$ are orthogonal local martingales under
$P^{\theta_0}$. 
The assumed boundedness further implies that they are martingales.
This implies that 
$\mathcal{E}_t := \exp\{L_t(\theta,\theta_0)\}$ is a
martingale under $P^{\theta_0}$.

Define $U^{\theta,f}$ by
\begin{equation*}
U^{\theta, f}_t =  f(Z_t) -f(Z_0) - \int_0^t \mathcal{L}^\theta 
f(Z_s)\mathrm{d}s
\end{equation*}
for $f \in C_b^{0,2}$.
If $\mathcal{E} U^{\theta,f}$ is a martingale under $P^{\theta_0}$, 
it means that 
$U^{\theta,f}$ is a martingale under the probability measure $P^\prime$
defined by
\begin{equation*}
 \frac{\mathrm{d}Q}{\mathrm{d}P^{\theta_0}} = \mathcal{E}_T.
\end{equation*}
If this is true for any $f$, $P^\prime$ is the solution of the martingale
problem with respect to the generator $\mathcal{L}^\theta$ 
and so, by Theorem~\ref{thm0}, $P^\prime = P^\theta$ and 
we get the result.
Therefore it only remains to show that 
$\mathcal{E} U^{\theta,f}$ is a martingale under $P^{\theta_0}$ for any
$f \in C_b^{0,2}$.
By It\^o's formula,
\begin{equation*}
 \mathcal{E}_\tau U^{\theta,f}_\tau = 
\int_0^\tau\mathcal{E}_t\mathrm{d}U^{\theta,f}_t +
\int_0^\tau U^{\theta,f}_t\mathrm{d}\mathcal{E}_t
+ \int_0^\tau \mathcal{E}_t\langle \mu^\theta-\mu^{\theta_0}, \nabla_y
f\rangle(Z_t)\mathrm{d}t + 
\sum_{0 < t \leq \tau} \Delta \mathcal{E}_t\Delta U^{\theta,f}_t
\end{equation*}
and
\begin{equation*}
\Delta \mathcal{E}_t = \mathcal{E}_{t-}\sum_{i\neq j}
\left(
\frac{q^{\theta}_{ji}(Y_t)}{q^{\theta_0}_{ji}(Y_t)} -1
\right)(N^{ji}_t-N^{ji}_{t-}),\ \ 
\Delta U^{\theta,f}_t = f(X_t,Y_t)-f(X_{t-},Y_t).
\end{equation*}
Since
\begin{equation*}
 \Delta \mathcal{E}_t \Delta U^{\theta,f}_t = 
\mathcal{E}_{t-}
\sum_{i\neq j}
\left(
\frac{q^{\theta}_{ji}(Y_t)}{q^{\theta_0}_{ji}(Y_t)} -1
\right)(N^{ji}_t-N^{ji}_{t-})(f(e_j,Y_t)-f(e_i,Y_t)),
\end{equation*}
we have
\begin{equation*}
\begin{split}
&\sum_{0 < t \leq \tau}  \Delta \mathcal{E}_t \Delta U^{\theta,f}_t \\
& =
\int_0^\tau
\mathcal{E}_{t-}
\sum_{i \neq j}
\left(
\frac{q^{\theta}_{ji}(Y_t)}{q^{\theta_0}_{ji}(Y_t)} -1
\right)(f(e_j,Y_t)-f(e_i,Y_t)) \mathrm{d}N^{ji}_t
 \\
& = \int_0^\tau
\mathcal{E}_{t-}
\sum_{i\neq j}
\left(
\frac{q^{\theta}_{ji}(Y_t)}{q^{\theta_0}_{ji}(Y_t)} -1
\right)(f(e_j,Y_t)-f(e_i,Y_t)) [\mathrm{d}N^{ji}_t -
 q^{\theta_0}_{ji}(Y_t)\langle X_t,e_i \rangle \mathrm{d}t] 
\\ & \hspace*{0.5cm} + 
\int_0^\tau
\mathcal{E}_t
\sum_{i \neq j}
\left(
q^{\theta}_{ji}(Y_t) - q^{\theta_0}_{ji}(Y_t)
\right)(f(e_j,Y_t)-f(e_i,Y_t))\langle X_t,e_i \rangle\mathrm{d}t.
\end{split}
\end{equation*}
Consequently, we have
\begin{equation*}
\begin{split}
  \mathcal{E}_\tau U^{\theta,f}_\tau = &
\int_0^\tau\mathcal{E}_t\mathrm{d}U^{\theta_0,f}_t +
\int_0^\tau U^{\theta,f}_t\mathrm{d}\mathcal{E}_t
\\ & \int_0^\tau
\mathcal{E}_{t-}
\sum_{i\neq j}
\left(
\frac{q^{\theta}_{ji}(Y_t)}{q^{\theta_0}_{ji}(Y_t)} -1
\right)(f(e_j,Y_t)-f(e_i,Y_t)) [\mathrm{d}N^{ji}_t -
 q^{\theta_0}_{ji}(Y_t)\langle X_t,e_i \rangle\mathrm{d}t],
\end{split}
\end{equation*}
which is a martingale under $P^{\theta_0}$ by (\ref{nji2}). \hfill////
\\

Theorem~\ref{thm2} provides a starting point of asymptotic theories.
Let us discuss briefly the large sample asymptotics $T\to \infty$ as an example.
Assume $\Theta \subset \mathbb{R}^n$ and the maps
$\theta \mapsto \mu^\theta$ and $\theta \mapsto Q^\theta =
[q^{\theta}_{ij}]$ 
are regular
enough. Denote by $\dot{\mu}^\theta$ and $\dot {q}_{ij}^\theta$ their
derivatives in $\theta$.
For $\theta_T = \theta_0 + uT^{-1/2}$, $u \in \mathbb{R}^n$,
by the Taylor expansion, we have
\begin{equation*}
\begin{split}
& \int_0^T\log
  \frac{q^{\theta_T}_{ji}(Y_t)}{q^{\theta_0}_{ji}(Y_t)} \mathrm{d}N^{ji}_t
- \int_0^T (q^{\theta_T}_{ji}(Y_t) - q^{\theta_0}_{ji}(Y_t))\langle
  X_t,e_i\rangle \mathrm{d}t
\\ &\approx \int_0^T
\log \frac{q^{\theta_T}_{ji}(Y_t)}{q^{\theta_0}_{ji}(Y_t)}
 [\mathrm{d}N^{ji}_t
-q^{\theta_0}_{ji}(Y_t)\langle
  X_t,e_i\rangle \mathrm{d}t] - \frac{1}{2}
\int_0^T \left|
\log \frac{q^{\theta_T}_{ji}(Y_t)}{q^{\theta_0}_{ji}(Y_t)}\right|^2 
q^{\theta_0}_{ji}(Y_t)\langle
  X_t,e_i\rangle \mathrm{d}t\\
&\approx
\frac{1}{\sqrt{T}}\int_0^T
 \frac{\langle u,\dot{q}^{\theta_0}_{ji}\rangle}{q^{\theta_0}_{ji}}(Y_t)
 [\mathrm{d}N^{ji}_t
-q^{\theta_0}_{ji}(Y_t)\langle
  X_t,e_i\rangle \mathrm{d}t] - \frac{1}{2T}
\int_0^T
 \frac{|\langle u,\dot{q}^{\theta_0}_{ji}\rangle|^2}{q^{\theta_0}_{ji}}(Y_t)
\langle
  X_t,e_i\rangle \mathrm{d}t
\end{split}
\end{equation*}
and
\begin{equation*}
 \begin{split}
&\frac{1}{\epsilon^2}\int_0^T \langle
  \mu^{\theta_T}(Z_t)-\mu^{\theta_0}(Z_t),\mathrm{d}Y_t\rangle
- \frac{1}{2\epsilon^2} \int_0^T (|\mu^{\theta_T}(Z_t)|^2 -
  |\mu^{\theta_0}(Z_t)|^2)\mathrm{d}t\\
& =
\frac{1}{\epsilon^2}\int_0^T \langle
  \mu^{\theta_T}(Z_t)-\mu^{\theta_0}(Z_t),\mathrm{d}Y_t-
\mu^{\theta_0}(Z_t)\mathrm{d}t\rangle
- \frac{1}{2\epsilon^2} \int_0^T |\mu^{\theta_T}(Z_t) -
  \mu^{\theta_0}(Z_t)|^2\mathrm{d}t\\
& \approx
\frac{1}{\epsilon\sqrt{T}}\int_0^T \Biggl\langle
  \langle u,\dot{\mu}^{\theta_0}(Z_t)\rangle,
\frac{\mathrm{d}Y_t-
\mu^{\theta_0}(Z_t)\mathrm{d}t}{\epsilon}\Biggr\rangle
- \frac{1}{2\epsilon^2T} \int_0^T | \langle u,
  \dot{\mu}^{\theta_0}(Z_t)\rangle|^2\mathrm{d}t.
 \end{split}
\end{equation*}
Therefore, under the ergodic property (see Yin and Zhu~\cite{YZ}):
\begin{equation*}
 \frac{1}{T}
\int_0^T
 \frac{|\langle u,\dot{q}^{\theta_0}_{ji}\rangle|^2}{q^{\theta_0}_{ji}}(Y_t)
\langle
  X_t,e_i\rangle \mathrm{d}t \to \int
\frac{|\langle u,\dot{q}^{\theta_0}_{ji}\rangle|^2}{q^{\theta_0}_{ji}}
(y)\langle x, e_i \rangle \pi(\mathrm{d}x,\mathrm{d}y),
\end{equation*}
\begin{equation*}
 \frac{1}{T} \int_0^T | \langle u,
  \dot{\mu}^{\theta_0}(Z_t)\rangle|^2\mathrm{d}t \to
\int |\langle u,\dot{\mu}^{\theta_0} \rangle|^2 (x,y)\pi(\mathrm{d}x,\mathrm{d}y),
\end{equation*}
the martingale central limit theorem implies the Local
Asymptotic Normality (LAN) property:
\begin{equation*}
 L_T(\theta_T,\theta_0) \to \langle u, N\rangle - \frac{1}{2}\langle u, I(\theta_0) u
  \rangle, \ \ N \sim \mathcal{N}(0,I(\theta_0)),
\end{equation*}
where
\begin{equation*}
 I(\theta_0) = \int \sum_{i \neq  j}
\frac{\dot{q}^{\theta_0}_{ji}(\dot{q}^{\theta_0}_{ji})^{\perp}}{q^{\theta_0}_{ji}}
(y)\langle x, e_i \rangle
+ \frac{1}{\epsilon^2}\dot{\mu}^{\theta_0} (\dot{\mu}^{\theta_0})^{\perp}(x,y)
 \pi(\mathrm{d}x,\mathrm{d}y)
\end{equation*}
and $\pi$ is the ergodic distribution. Under the non-degeneracy of
$I(\theta_0)$, the LAN implies the asymptotic efficiency of a Maximum
Likelihood Estimator (MLE), with asymptotic variance $I(\theta_0)^{-1}$;
see Ibragimov and Has'minskii~\cite{IH} for the detail.

\section{A finite-dimensional filter}
Here we extend the filtering theory of hidden Markov models developed in
Chapter 8 of Elliott et al.~\cite{El} to the SHS
\begin{equation*}
 \begin{split}
&  \mathrm{d}X_t = Q(Y_t)X_t\mathrm{d}t + \mathrm{d}V_t,\\
& \mathrm{d}Y_t = \mu(X_t,Y_t)\mathrm{d}t + \epsilon \mathrm{d}W_t,
 \end{split}
\end{equation*}
where $V$ is an $\{\mathcal{F}_t\}$-martingale (recall Corollary~\ref{cor1}).
Denoting $\langle e_i,\mu(e_j,y)\rangle = c_{ij}(y)$, $C(y) = [c_{ij}(y)]$, 
we can write $\mu(X_s,Y_s) = C(Y_s)X_s$.
In this section we assume we observe only a continuous sample path 
$Y$ on a time interval $[0,T]$
while $X$ is hidden.
The system is a hidden Markov model in \cite{El} when both $Q$ and $\mu$
do not depend on $Y$.
By the dependence, $V$ is not independent of $W$ and so,
it is beyond the framework of \cite{El}.
We however show in this and the next sections 
that the results in Chapter 8 of \cite{El} remain valid.
Namely, a finite-dimensional  filter and the
 EM algorithm  can be constructed for
the SHS. A key for this is Lemma~\ref{lem2}.\\

Denote by $\mathcal{F}^Y$ the natural filtration of $Y$.
The filtering problem is to infer $X$ from the observation of $Y$,
that is, to compute $E[X_t|\mathcal{F}^Y_t]$.
The smoothing problem is to compute 
$E[X_t|\mathcal{F}^Y_T]$ for $t \leq T$.
Denote $E^0_t[H] = E^0[H|\mathcal{F}^Y_t]$ for a given integrable random
variable $H$, where $E^0$ is the expectation under $P^0$ in Section~2.
For a given process $H$,
the Bayes formula gives
\begin{equation}\label{Bayes}
 E[H_t | \mathcal{F}^Y_t] = \frac{E^0_t[\Lambda_t H_t]}
{E^0_t[\Lambda_t]},
\end{equation}
where $\Lambda$ is defined by (\ref{Lambda}).
Theorem~\ref{thm42} below shows that both the numerator and denominator
of (\ref{Bayes}) for $H_t = X_t$ can be computed by solving a linear equation.

\begin{thm}\label{thm3}
Under the same conditions of Theorem~\ref{thm1},
if $H$ is  of the form
\begin{equation}\label{Hform}
 \mathrm{d}H_t = \alpha_t \mathrm{d}t + \langle \beta_t, \mathrm{d}X_t
  \rangle + \langle \delta_t,  \mathrm{d}Y_t \rangle,
\end{equation}
where $\alpha, \beta, \delta$ are bounded predictable processes, 
then
\begin{equation}\label{filt}
\begin{split}
 E^0_t[\Lambda_tH_t]  = &  H_0 + \frac{1}{\epsilon^2}
  \int_0^t \langle C(Y_s)E^0_s[\Lambda_sH_sX_s] + \epsilon^2 
E^0_s[\Lambda_s\delta_s], \mathrm{d}Y_s\rangle
 \\ & + \int_0^t E^0_s[\Lambda_s\alpha_s] +
E^0_s[\Lambda_s\langle \beta_s, Q(Y_s)X_{s-}\rangle] +
 E^0_s[\Lambda_s\langle\delta_s, C(Y_s)X_s\rangle]\mathrm{d}s.
\end{split}
\end{equation}
\end{thm}
{\it Proof: }
It\^o's formula gives
\begin{equation*}
 \Lambda_t H_t = H_0 + \int_0^t H_{s-} \mathrm{d}\Lambda_s 
+ \int_0^t \Lambda_s \mathrm{d}H_s + \int_0^t \Lambda_s \langle \delta_s,
\mu(Z_s)\rangle \mathrm{d}s.
\end{equation*}
Take the conditional expectation under $P^0$ given $\mathcal{F}^Y_t$ to 
get (\ref{filt}).
Here, we have used the fact that $Y/\epsilon$ is a $d$ dimensional Brownian
motion under $P^0$ as well as Lemma~\ref{lem2}.
\hfill////

\begin{thm} \label{thm42}
 Under the same conditions of Theorem~\ref{thm1},
 \begin{equation}\label{f1}
 E^0_t[\Lambda_tX_t] = 
 X_0 + \frac{1}{\epsilon^2}\int_0^t 
\mathrm{diag}(E^0_s[\Lambda_sX_s])C(Y_s)^{\perp}\mathrm{d}Y_s
+ \int_0^t Q(Y_s) E^0_s[\Lambda_s X_s]  \mathrm{d}s
 \end{equation}
and
\begin{equation}\label{f2}
 E^0_t[\Lambda_t] =
\langle \mathbf{1}, E^0_t[\Lambda_tX_t]\rangle
= 1 + \frac{1}{\epsilon^2}\int_0^t \langle
  C(Y_s)E^0_s[\Lambda_sX_s] ,\mathrm{d}Y_s \rangle.
\end{equation}
\end{thm}
{\it Proof: } 
For each $i = 1,\dots,k$, let $H_t = \langle e_i,X_t\rangle$ in (\ref{filt})
to get
 \begin{equation*}
\begin{split}
&  \langle e_i, E^0_t[\Lambda_tX_t]\rangle \\ &=
 \langle e_i,X_0 \rangle
+ \frac{1}{\epsilon^2}\int_0^t \langle e_i, E^0_s[\Lambda_sX_s]\rangle
 \langle C(Y_s)e_i,  \mathrm{d}Y_s\rangle
+ \int_0^t \langle e_i, Q(Y_s) E^0_s[\Lambda_s X_s] \rangle \mathrm{d}s,
\end{split}
 \end{equation*}
which is equivalent to (\ref{f1}). 
Here we have used that $\langle e_i, X_s\rangle X_s = 
\langle e_i, X_s\rangle e_i$.
The first identity of (\ref{f2}) is by $\langle \mathbf{1},X_t\rangle =
1$.
To get the second identity, use the first and (\ref{f1}), or alternatively,
let $H_t = 1$ in (\ref{filt}).
\hfill////
\\

The smoothing problem is also solved as follows.

\begin{thm}
Under the same conditions of Theorem~\ref{thm1}, for $t \geq \tau$,
\begin{equation*}
 E^0_t[\Lambda_tX_\tau] =  E^0_\tau[\Lambda_\tau X_\tau] 
 + \frac{1}{\epsilon^2}\int_{\tau}^t
\mathrm{diag}(E^0_s[\Lambda_sX_\tau])C(Y_s)^{\perp}\mathrm{d}Y_s.
\end{equation*}
\end{thm}
{\it Proof: } Let $H_t = \langle e_i, X_{t \wedge \tau} \rangle$ in
(\ref{filt}) to get
 \begin{equation*}
\langle e_i, E^0_t[\Lambda_tX_\tau]\rangle  = 
 \langle e_i,E^0_\tau[\Lambda_\tau X_\tau]\rangle
+ \frac{1}{\epsilon^2}\int_{\tau}^t \langle e_i, E^0_s[\Lambda_sX_\tau]\rangle
 \langle C(Y_s)e_i,  \mathrm{d}Y_s\rangle
 \end{equation*}
for $t\geq \tau$. 
 \hfill////\\

\section{The likelihood under partial observations}
Here we consider again the parametric family $\{P^\theta\}$ introduced
in Section~3. We assume that a continuous sample path $Y$ is observed
on a time interval $[0,T]$ while $X$ is hidden. 
Under the same assumptions as in Theorem~\ref{thm2},
the law of $Y$ under $P^\theta$ is equivalent to that under
$P^{\theta_0}$
and the log likelihood function is given by
\begin{equation}\label{LY}
 L^Y(\theta,\theta_0) = \log E^{\theta_0}\left[
\frac{\mathrm{d}P^{\theta}}{\mathrm{d}P^{\theta_0}} \Big| \mathcal{F}^Y_T
\right].
\end{equation}
A MLE is therefore given by
\begin{equation*}
\hat{\theta} =  \mathrm{argmax}_{\theta \in \Theta} L^Y(\theta,\theta_0).
\end{equation*} 
Note that $\hat{\theta}$ does not depend on the choice of $\theta_0$
because by the Bayes formula,
\begin{equation}\label{choice}
\begin{split}
 L^Y(\theta,\theta_0)
& = \log  E^{\theta_1}
\left[ \frac{\mathrm{d}P^{\theta_0}}{\mathrm{d}P^{\theta_1}}
\frac{\mathrm{d}P^{\theta}}{\mathrm{d}P^{\theta_0}} \Big| \mathcal{F}^Y_T
\right]
- \log  E^{\theta_1}
\left[ \frac{\mathrm{d}P^{\theta_0}}{\mathrm{d}P^{\theta_1}}
\Big| \mathcal{F}^Y_T
\right] \\
& = L^Y(\theta,\theta_1) - L^Y(\theta_0,\theta_1) 
\end{split}
\end{equation}
for any $\theta_1 \in \Theta$. 

Here we show that the likelihood function can be obtained as a solution of a
linear filtering equation. Let
\begin{equation*}
H_t = \frac{\mathrm{d} P^{\theta}}{\mathrm{d}P^{\theta_0}}
 \Big|_{\mathcal{F}_t} = \exp(L_t(\theta,\theta_0)).
\end{equation*}
By Theorem~\ref{thm2}, we have
\begin{equation*}
\begin{split}
 \mathrm{d}H_t = H_{t-}\bigl\{ &
\langle R^{\theta,\theta_0}(Y_t)X_{t-}, \mathrm{d}X_t -
Q^{\theta_0}(Y_t)X_t\mathrm{d}t\rangle \\ & + \frac{1}{\epsilon^2}
\langle(C^\theta(Y_t)-C^{\theta_0}(Y_t))X_t, \mathrm{d}Y_t -
C^{\theta_0}(Y_t)X_t \mathrm{d}t \rangle
\bigr\},
\end{split}
\end{equation*}
where
\begin{equation*}
 R^{\theta,\theta_0}(Y_t) = \left[r^{\theta,\theta_0}_{ij}(Y_t)\right],
\ \ r^{\theta,\theta_0}_{ij}(Y_t) = \begin{cases}
\frac{q_{ij}^\theta}{q_{ij}^{\theta_0}}(Y_t)-1, \ \ & i \neq j
 \\
	0, \ \ & i=j.			     
				    \end{cases}
\end{equation*}
The likelihood function under the partial observation of $\{Y_t\}_{t \in
[0,T]}$ is given by
\begin{equation}\label{Lh}
\exp(L^Y(\theta,\theta_0)) = 
 E^{\theta_0}\left[H_T | \mathcal{F}^Y_T\right]
= \frac{E^0_T[\Lambda_T H_T]}{E^0_T[\Lambda_T]}.
\end{equation}
Both the numerator and denominator are computed as follows.
\begin{thm}\label{thm61}
Under the same conditions of Theorem~\ref{thm2},  for $t \in [0,T]$,
\begin{equation}\label{new1}
 E^0_t[\Lambda_tH_t] =\langle \mathbf{1},E^0_t[\Lambda_tH_tX_t]\rangle  =
  1 + \frac{1}{\epsilon^2}
  \int_0^t \langle C^{\theta}(Y_s)E^0_s[\Lambda_sH_sX_s], 
\mathrm{d}Y_s\rangle
\end{equation}
and 
 \begin{equation}\label{new2}
\begin{split}
&E^0_t[\Lambda_tH_tX_t] \\&= 
 X_0 + \frac{1}{\epsilon^2}\int_0^t 
\mathrm{diag}(E^0_s[\Lambda_sH_sX_s])C^\theta(Y_s)^{\perp}\mathrm{d}Y_s
+ \int_0^t Q^\theta(Y_s) E^0_s[\Lambda_s H_sX_s]  \mathrm{d}s.
\end{split}
 \end{equation}
\end{thm}
{\it Proof: }
Since
\begin{equation*}
\begin{split}
 \Delta H_t \langle e_i, \Delta X_t\rangle &=
H_{t-}\sum_{j,l=1}^k
\langle X_{t-},e_l\rangle r^{\theta,\theta_0}_{jl}(Y_t)\langle
e_j,\Delta X_t \rangle \langle e_i, \Delta X_t\rangle
 \\
& = H_{t-}\left\{
-\sum_{j=1}^k\langle X_{t-},e_i\rangle r^{\theta,\theta_0}_{ji}(Y_t)\langle
e_j,\Delta X_t \rangle + \sum_{l=1}^k 
\langle X_{t-},e_l\rangle r^{\theta,\theta_0}_{il}(Y_t)\langle e_i, \Delta X_t\rangle
\right\}\\
&= H_{t-}\left\{
-\langle X_{t-},e_i\rangle \langle R^{\theta,\theta_0}X_{t-}
,\Delta X_t \rangle + 
\langle e_i, R^{\theta,\theta_0}(Y_t)X_{t-}\rangle\langle e_i, \Delta X_t\rangle
\right\},
\end{split}
\end{equation*}
applying Theorem~\ref{thm3} with
\begin{equation*}
\begin{split}
 &\alpha_t = -\langle e_i,X_t\rangle H_t\left\{
\langle R^{\theta,\theta_0}(Y_t)X_t,Q^{\theta_0}(Y_t)X_t\rangle
+ \frac{1}{\epsilon^2}\langle 
(C^\theta(Y_t)-C^{\theta_0}(Y_t))X_t, C^{\theta_0}(Y_t)X_t \rangle
\right\},\\
& \beta_t= 
H_{t-}e_i  + H_{t-} \langle e_i, R^{\theta,\theta_0}(Y_t)X_{t-}\rangle e_i,\\
& \delta_t= \frac{1}{\epsilon^2} \langle e_i,X_t\rangle
H_t(C^\theta(Y_t)-C^{\theta_0}(Y_t))X_t,
\end{split}
\end{equation*}
we obtain, using $\langle e_i,X_s\rangle X_s = \langle e_i,X_s\rangle
e_i$ again and again,
\begin{equation*}
\begin{split}
 \langle e_i, &E^0_t[\Lambda_tH_tX_t] \rangle\\
=& \langle e_i,X_0\rangle + 
\frac{1}{\epsilon^2} \int_0^t \langle e_i, E^0_s[\Lambda_sH_sX_s]\rangle 
\langle C^{\theta}(Y_s)e_i, \mathrm{d}Y_s\rangle
 \\ & + \int_0^t E^0_s[\Lambda_s\alpha_s] +
E^0_s[\Lambda_s\langle \beta_s, Q^{\theta_0}(Y_s)X_{s-}\rangle] +
 E^0_s[\Lambda_s\langle\delta_s, C^{\theta_0}(Y_s)X_s\rangle]\mathrm{d}s\\
=& \langle e_i,X_0\rangle + 
\frac{1}{\epsilon^2} \int_0^t \langle e_i, E^0_s[\Lambda_sH_sX_s]\rangle 
\langle C^{\theta}(Y_s)e_i, \mathrm{d}Y_s\rangle
 \\ & - \int_0^t 
\langle R^{\theta,\theta_0}(Y_s)e_i,Q^{\theta_0}(Y_s)e_i \rangle
\langle e_i, E^0_s[\Lambda_sH_sX_s] \rangle \mathrm{d}s \\
& + \int_0^t E^0_s[\Lambda_sH_s(\langle e_i,Q^{\theta_0}(Y_s)X_s\rangle
+ \langle e_i,R^{\theta,\theta_0}(Y_s)X_s\rangle\langle
e_i,Q^{\theta_0}(Y_s)X_s\rangle)]\mathrm{d}s.
\end{split}
\end{equation*}
Since
\begin{equation*}
 \langle R^{\theta,\theta_0}(Y_s)e_i,Q^{\theta_0}(Y_s)e_i \rangle
= \sum_{j=1} 
r^{\theta,\theta_0}_{ji}(Y_s)q^{\theta_0}_{ji}(Y_s) =
\sum_{j \neq i} (q^{\theta}_{ji} -
q^{\theta_0}_{ji})(Y_s)
= - (q^{\theta}_{ii} -
q^{\theta_0}_{ii})(Y_s)
\end{equation*}
and
\begin{equation*}
\begin{split}
 \langle e_i &,R^{\theta,\theta_0}(Y_s)X_s\rangle\langle
e_i,Q^{\theta_0}(Y_s)X_s\rangle\\
&= \sum_{j=1}^k\langle e_i,R^{\theta,\theta_0}(Y_s)e_j\rangle\langle
e_i,Q^{\theta_0}(Y_s)e_j\rangle \langle e_j,X_s\rangle
\\& =  
\sum_{j\neq i}\langle e_i,
(Q^\theta(Y_t) -Q^{\theta_0}(Y_s))e_j\rangle \langle e_j,X_s\rangle
\\
&= \langle e_i,
(Q^\theta(Y_t) -Q^{\theta_0}(Y_s))X_s\rangle 
- (q_{ii}^\theta(Y_t) -q_{ii}^{\theta_0}(Y_s))\langle e_i,X_s \rangle ,
\end{split}
\end{equation*}
we conclude the linear equation 
\begin{equation*}
\begin{split}
 \langle e_i, E^0_t[\Lambda_tH_tX_t] \rangle = \langle e_i,X_0\rangle  
\ +  &  \ \frac{1}{\epsilon^2} \int_0^t \langle e_i, E^0_s[\Lambda_sH_sX_s]\rangle 
\langle C^{\theta}(Y_s)e_i, \mathrm{d}Y_s\rangle
\\
& + \int_0^t \langle e_i,Q^{\theta}(Y_s)E^0_s[\Lambda_sH_sX_s]\rangle
\mathrm{d}s,
\end{split}
\end{equation*}
which is equivalent to  (\ref{new2}). Then,
(\ref{new1}) follows from $\langle \mathbf{1},X_t\rangle = 1$.
 \hfill////

\begin{cor}
Let $L^Y(\theta,\theta_0)$ be
the log likelihood function defined by (\ref{LY}), and
\begin{equation*}
 \hat{\mu}^\theta_t = E^\theta[\mu^\theta(Z_t)|\mathcal{F}^Y_t].
\end{equation*}
Then,
\begin{equation*}
\begin{split}
L^Y(\theta,\theta_0) & = 
\frac{1}{\epsilon^2}
\int_0^T\langle  
\hat{\mu}^{\theta}_t - \hat{\mu}^{\theta_0}_t , \mathrm{d}Y_t \rangle -
\frac{1}{2\epsilon^2}\int_0^T  
|\hat{\mu}^{\theta}_t|^2 - |\hat{\mu}^{\theta_0}_t|^2
 \mathrm{d}t \\
& = 
\frac{1}{\epsilon^2}
\int_0^T\langle  
\hat{\mu}^{\theta}_t - \hat{\mu}^{\theta_0}_t , \mathrm{d}Y_t - \hat{\mu}^{\theta_0}_t\mathrm{d}t\rangle -
\frac{1}{2\epsilon^2}\int_0^T  
|\hat{\mu}^{\theta}_t-\hat{\mu}^{\theta_0}_t|^2
 \mathrm{d}t.
\end{split}
\end{equation*}
\end{cor}
{\it Proof: } First we should note that this representation is valid in
 a more general framework than the stochastic hybrid system. Indeed, by Theorem
~7.12 of Liptser and Shiryaev~\cite{LS1},
\begin{equation*}
\hat{W}_t = \frac{1}{\epsilon}\left\{ Y_t - \int_0^t
			       \hat{\mu}^{\theta}_s\mathrm{d}s \right\}
\end{equation*}
is an $\left\{\mathcal{F}^Y_t\right\}$-Brownian motion under
$P^{\theta}$.
Since
\begin{equation}\label{diffmod}
 \mathrm{d}Y_t = \hat{\mu}^{\theta}_t \mathrm{d}t + \epsilon \mathrm{d}\hat{W}_t,
\end{equation}
the result follows from the Girsanov-Maruyama theorem.
Here we derive it as a corollary of Theorem~\ref{thm61}.
By (\ref{Lh}), we have
\begin{equation*}
 L^Y(\theta,\theta_0) = \log E^0_T[\Lambda_TH_T] - \log E^0_T[\Lambda_T],
\end{equation*}
and by It\^o's formula and (\ref{new1}),
\begin{equation*}
\begin{split}
 \mathrm{d}\log E^0_t[\Lambda_tH_t]
&= \frac{1}{\epsilon^2}\Biggl\langle
C^\theta(Y) \frac{E^0_t[\Lambda_tH_tX_t]}{E^0_t[\Lambda_tH_t]},
\mathrm{d}Y_t \Biggr\rangle
- \frac{1}{2\epsilon^2}
\left|C^\theta(Y)
 \frac{E^0_t[\Lambda_tH_tX_t]}{E^0_t[\Lambda_tH_t]}\right|^2 \mathrm{d}t
\\
&=  
\frac{1}{\epsilon^2}\langle
C^\theta(Y)E^{\theta}[X_t|\mathcal{F}^Y_t]
\mathrm{d}Y_t 
- \frac{1}{2\epsilon^2}
|C^\theta(Y)E^{\theta}[X_t|\mathcal{F}^Y_t]|^2 \mathrm{d}t.
\end{split}
\end{equation*}
Similarly, we obtain
\begin{equation*}
 \mathrm{d}\log E^0_t[\Lambda_t]
=  
\frac{1}{\epsilon^2}\langle
C^{\theta_0}(Y)E^{\theta_0}[X_t|\mathcal{F}^Y_t]
\mathrm{d}Y_t 
- \frac{1}{2\epsilon^2}
|C^{\theta_0}(Y)E^{\theta_0}[X_t|\mathcal{F}^Y_t]|^2 \mathrm{d}t
\end{equation*}
and conclude the result.\hfill////
\\

By (\ref{new1}), (\ref{new2}), (\ref{diffmod}),  and
\begin{equation*}
 \hat{\mu}^\theta_t = C^\theta(Y_t)E^{\theta}[X_t|\mathcal{F}^Y_t]
= C^\theta(Y_t) \frac{\bar{X}^\theta_t}{\langle \mathbf{1},\bar{X}^\theta_t\rangle},\
\ 
\bar{X}^\theta_t = E^0_t[\Lambda_tH_tX_t],
\end{equation*}
we observe that our model is translated to a path-dependent diffusion model
\begin{equation*}
 \begin{split}
  &\mathrm{d}Y_t = C^\theta(Y_t) \frac{\bar{X}^\theta_t}{\langle
  \mathbf{1},\bar{X}^\theta_t\rangle}\mathrm{d}t + \epsilon
  \mathrm{d}\hat{W}_t, \\
& \mathrm{d}\bar{X}^\theta_t =  Q^\theta(Y_t)\bar{X}^\theta_t
  \mathrm{d}t 
+  \frac{1}{\epsilon^2}
\mathrm{diag}(\bar{X}^\theta_t)C^\theta(Y_t)^\perp \mathrm{d}Y_t.
 \end{split}
\end{equation*}
This representation would serve as a starting point of a theory of
asymptotic statistics, which however remains for future research.
In particular, specifications of 
an appropriate ergodicity condition,
an identifiability condition, and the asymptotic variance of a MLE  remain open. 
Nevertheless, having the complete observation counterpart in our mind
(see Section~3),
it would be natural to expect a  MLE
to work reasonably well also under the partial observation. In the sequel, we
focus on a computational issue of a MLE.

An apparent computational problem to obtain a MLE:
\begin{equation*}
\hat{\theta} =  \mathrm{argmax}_{\theta \in \Theta} L^Y(\theta,\theta_0)
= \mathrm{argmax}_{\theta \in \Theta} \langle
\mathbf{1},\bar{X}^\theta_T \rangle
\end{equation*} 
 is that the associated numerical maximization
is time-consuming, due to that
the evaluation of the objective function requires to solve  a different
filtering equation for each $\theta$.
For some cases however an efficient algorithm is available as will
be seen in the next section.
\\

\section{The EM algorithm}

The EM algorithm is a well-known technique which is often useful when a
likelihood function is not easily computable.
Let us recall the idea of the EM algorithm.
Let
\begin{equation*}
 Q(\theta^\ast,\theta) = 
  E^{\theta}\left[ \log
\frac{\mathrm{d}P^{\theta^\ast}}{\mathrm{d}P^{\theta}} \Big| \mathcal{F}^Y_T
\right].
\end{equation*}
By 
Jensen's inequality and (\ref{choice}), 
\begin{equation}\label{em}
Q(\theta^\ast,\theta)
\leq
 \log E^{\theta}\left[
\frac{\mathrm{d}P^{\theta^\ast}}{\mathrm{d}P^{\theta}} \Big| \mathcal{F}^Y_T
\right] = L^Y(\theta^\ast,\theta) = L^Y(\theta^\ast,\theta_0) - L^Y(\theta,\theta_0),
\end{equation}
which means that the sequence defined by
\begin{equation*}
 \theta_{n+1} = \mathrm{argmax}_{\theta \in \Theta} Q(\theta,\theta_n)
\end{equation*}
makes $L^Y(\theta_n,\theta_0)$ increasing.
Under an appropriate condition the sequence $\{\theta_n\}$ converges to 
a MLE  $\hat{\theta}$, for which we refer to
Wu~\cite{Wu}. See also Dembo and Zeitouni~\cite{DZ}.
As easily seen from (\ref{em}), a MLE is a stationary point of the EM algorithm in
general.
However in general, there might be many stationary points other than a
MLE and so, as in many other optimization problems, 
one would need to examine results with various initial parameters $\theta_0$.

The EM algorithm is particularly 
effective when the statistical model in question is
an exponential family under complete observations.
Let us consider the case where
 $Q^\theta(y) = [q^\theta_{ji}(y)]$ and $C^\theta(y)$ are of the form
\begin{equation*}
 q^{\theta}_{ji}(y) = \varphi_{ji}(\theta)q^0_{ji}(y),\ \ 
\mu^\theta(z) = \sum_{l=1}^L \psi_l(\theta) \mu^l(z)
\end{equation*}
for some functions $\varphi_{ji}, \psi_l : \Theta \to \mathbb{R}$
and $k\times k$ Q-matrix $Q^0(y) = [q^0_{ji}(y)]$ and 
$\mu^l : S \mapsto \mathbb{R}^d$, $l=1,\dots L$.
The log likelihood under the complete observation (see
Theorem~\ref{thm2}) is then
\begin{equation*}
\begin{split}
 L_T(\theta,\theta_0) = & \sum_{i\neq j} \left\{
N^{ji}_T
\log \frac{\varphi_{ji}(\theta)}{\varphi_{ji}(\theta_0)}
- \left(\frac{\varphi_{ji}(\theta)}{\varphi_{ji}(\theta_0)}-1\right)
\int_0^T q_{ji}^{\theta_0}(Y_t)\langle e_i,X_t \rangle \mathrm{d}t
\right\} \\
&+\frac{1}{\epsilon^2}\sum_{l=1}^L
(\psi_l(\theta)-\psi_l(\theta_0))\int_0^T \langle \mu^l(Z_t),
\mathrm{d}Y_t - \mu^{\theta_0}(Z_t)\mathrm{d}t\rangle \\ 
& - \frac{1}{2\epsilon^2}\sum_{l,m=1}^L
(\psi_l(\theta)-\psi_l(\theta_0))(\psi_m(\theta)-\psi_m(\theta_0))
\int_0^T \langle \mu^l(Z_t),\mu^m(Z_t)\rangle \mathrm{d}t.
\end{split}
\end{equation*}
The statistical model under the complete observation 
is then an exponential family with sufficient statistics
\begin{equation*}
 N^{ji}_T, \ \ 
\int_0^T q_{ji}^{\theta_0}(Y_t)\langle e_i,X_t \rangle \mathrm{d}t,
\ \ 
\int_0^T \langle \mu^l(Z_t),
\mathrm{d}Y_t - \mu^{\theta_0}(Z_t)\mathrm{d}t\rangle, \ \ 
\int_0^T \langle \mu^l(Z_t),\mu^m(Z_t)\rangle \mathrm{d}t.
\end{equation*}
Under the partial observation, in terms of the EM algorithm,
the Expectation step amounts to filtering 
 the sufficient statistics under $P^{\theta_0}$:
\begin{equation*}
\begin{split}
 Q(\theta&,\theta_0) \\= & \sum_{i\neq j} \left\{
E^{\theta_0}[N^{ji}_T|\mathcal{F}^Y_T]
\log \frac{\varphi_{ji}(\theta)}{\varphi_{ji}(\theta_0)}
- \left(\frac{\varphi_{ji}(\theta)}{\varphi_{ji}(\theta_0)}-1\right)
E^{\theta_0}[\int_0^T q_{ji}^{\theta_0}(Y_t)\langle e_i,X_t \rangle \mathrm{d}t|\mathcal{F}^Y_T]
\right\} \\
&+\frac{1}{\epsilon^2}\sum_{l=1}^L
(\psi_l(\theta)-\psi_l(\theta_0))
E^{\theta_0}[\int_0^T \langle \mu^l(Z_t),\mathrm{d}Y_t - \mu^{\theta_0}(Z_t)\mathrm{d}t\rangle |\mathcal{F}^Y_T]
\\
& - \frac{1}{2\epsilon^2}\sum_{l,m=1}^L
(\psi_l(\theta)-\psi_l(\theta_0))(\psi_m(\theta)-\psi_m(\theta_0))
E^{\theta_0}[\int_0^T \langle \mu^l(Z_t),\mu^m(Z_t)\rangle \mathrm{d}t|\mathcal{F}^Y_T],
\end{split}
\end{equation*}
and the Maximization
step is the same exercise as finding a MLE for an exponential family.
The filtering equations to solve are as follows.

\begin{thm} 
Under the same condition of Theorem~\ref{thm2},
for each $i,j = 1,\dots,k$ and $l,m = 1,\dots,L$,
\begin{enumerate}
 \item 
\begin{equation*}
 E^{\theta_0}[N^{ji}_T | \mathcal{F}^Y_T] = \frac{E^0[\Lambda_T
  N^{ji}_T]}{E^0[\Lambda_T]}, \ \ 
E^0_T[\Lambda_TN^{ji}_T] = \langle \mathbf{1}, E^0_T[\Lambda_TN^{ji}_TX_T]
\rangle
\end{equation*}
and $F_t := E^0_t[\Lambda_tN^{ji}_tX_t] $ is the solution of
\begin{equation}\label{new3}
\begin{split}
F_t
= & \frac{1}{\epsilon^2}\int_0^t \mathrm{diag}(F_s)
C^{\theta_0}(Y_s)^\perp \mathrm{d}Y_s +  \int_0^t Q^{\theta_0}(Y_s)F_s \mathrm{d}s \\
& + e_j\int_0^t  q^{\theta_0}_{ji}(Y_s) \langle e_i,E^0_s[\Lambda_sX_s] \rangle
\mathrm{d}s.
\end{split}
\end{equation}
\item
\begin{equation*}
 E^{\theta_0}[
\int_0^T q_{ji}^{\theta_0}(Y_t)\langle e_i,X_t \rangle \mathrm{d}t
 | \mathcal{F}^Y_T] = \frac{1}{E^0_T[\Lambda_T]} {E^0_T[\Lambda_T
 \int_0^T q_{ji}^{\theta_0}(Y_t)\langle e_i,X_t \rangle \mathrm{d}t]},
\end{equation*}
\begin{equation*}
E^0_T[\Lambda_T
 \int_0^T q_{ji}^{\theta_0}(Y_t)\langle e_i,X_t \rangle \mathrm{d}t]
 = \langle \mathbf{1},
E^0_T[\Lambda_TX_T
 \int_0^T q_{ji}^{\theta_0}(Y_t)\langle e_i,X_t \rangle \mathrm{d}t]
\rangle
\end{equation*}
and 
\begin{equation*}
 F_t := E^0_t[\Lambda_tX_t
 \int_0^t q_{ji}^{\theta_0}(Y_s)\langle e_i,X_s \rangle \mathrm{d}s
]
\end{equation*}
is the solution of
\begin{equation*}
\begin{split}
 F_t = &\frac{1}{\epsilon^2}\int_0^t
  \mathrm{diag}(F_s)C^{\theta_0}(Y_s)^\perp \mathrm{d}Y_s 
+ \int_0^t Q^{\theta_0}(Y_s)F_s\mathrm{d}s \\
& +  e_i\int_0^t  q^{\theta_0}_{ji}(Y_s) \langle e_i,E^0_s[\Lambda_sX_s] \rangle
\mathrm{d}s.
\end{split}
\end{equation*}
\item
\begin{equation*}
\begin{split}
& E^{\theta_0}[
\int_0^T \langle \mu^l(Z_t),
\mathrm{d}Y_t - \mu^{\theta_0}(Z_t)\mathrm{d}t\rangle
 | \mathcal{F}^Y_T] \\ &= \frac{1}{E^0_T[\Lambda_T]} {E^0_T[\Lambda_T
 \int_0^T \langle \mu^l(Z_t),
\mathrm{d}Y_t - \mu^{\theta_0}(Z_t)\mathrm{d}t\rangle]},
\end{split}
\end{equation*}
\begin{equation*}
\begin{split}
&E^0_T[\Lambda_T
 \int_0^T \langle \mu^l(Z_t),
\mathrm{d}Y_t - \mu^{\theta_0}(Z_t)\mathrm{d}t\rangle]
\\ & = \langle \mathbf{1},
E^0_T[\Lambda_TX_T
 \int_0^T \langle \mu^l(Z_t),
\mathrm{d}Y_t - \mu^{\theta_0}(Z_t)\mathrm{d}t\rangle]
\rangle
\end{split}
\end{equation*}
and 
\begin{equation*}
 F_t := 
E^0_t[\Lambda_tX_t
 \int_0^t \langle \mu^l(Z_s),
\mathrm{d}Y_s - \mu^{\theta_0}(Z_s)\mathrm{d}s\rangle]
\end{equation*}
is the solution of
\begin{equation*}
\begin{split}
 F_t = &\frac{1}{\epsilon^2}\int_0^t
  \mathrm{diag}(F_s)C^{\theta_0}(Y_s)^\perp \mathrm{d}Y_s 
+ \int_0^t Q^{\theta_0}(Y_s)F_s\mathrm{d}s \\
& +  \int_0^t  \mathrm{diag}(E^0_s[\Lambda_sX_s])C^l(Y_s)^\perp \mathrm{d}Y_s,
\end{split}
\end{equation*}
where $C^l(y) = [c^l_{ij}]$, $c^l_{ij}(y) = \langle e_i,
     \mu^l(y,e_j))\rangle$.
\item
\begin{equation*}
E^{\theta_0}[
\int_0^T \langle \mu^l(Z_t),\mu^m(Z_t)\rangle \mathrm{d}t
 | \mathcal{F}^Y_T] = \frac{1}{E^0_T[\Lambda_T]} {E^0_T[\Lambda_T
 \int_0^T \langle \mu^l(Z_t),\mu^m(Z_t)\rangle \mathrm{d}t]},
\end{equation*}
\begin{equation*}
 E^0_T[\Lambda_T
 \int_0^T \langle \mu^l(Z_t),\mu^m(Z_t)\rangle \mathrm{d}t] = \langle \mathbf{1},
E^0_T[\Lambda_TX_T
 \int_0^T \langle \mu^l(Z_t),\mu^m(Z_t)\rangle \mathrm{d}t]
\rangle
\end{equation*}
and 
\begin{equation*}
 F_t := 
E^0_t[\Lambda_tX_t
\int_0^t \langle \mu^l(Z_s),\mu^m(Z_s)\rangle \mathrm{d}s]
\end{equation*}
is the solution of
\begin{equation*}
\begin{split}
 F_t = &\frac{1}{\epsilon^2}\int_0^t
  \mathrm{diag}(F_s)C^{\theta_0}(Y_s)^\perp \mathrm{d}Y_s 
+ \int_0^t Q^{\theta_0}(Y_s)F_s\mathrm{d}s \\
& +  \int_0^t  \mathrm{diag}(E^0_s[\Lambda_sX_s])D^{l,m}(Y_s)\mathrm{d}s,
\end{split}
\end{equation*}
where $D^{l,m}(y) =
\mathrm{diag}(\langle \mu^l(e_1,y), \mu^m(e_1,y)\rangle, 
\dots, \langle \mu^l(e_k,y), \mu^m(e_k,y)\rangle)$.
\end{enumerate} 
\end{thm}
{\it Proof: }
For 1., let $H_t = \langle e_a,X_t\rangle N^{ji}_t$ and apply
Theorem~\ref{thm3} with 
$\alpha = \delta = 0$ and $\beta_t = \delta_{aj} \langle
e_i,X_{t-}\rangle e_j + N^{ji}_{t-}e_a$ to obtain
\begin{equation*}
\begin{split}
 \langle e_a,E^0_t[&\Lambda_tN^{ji}_tX_t]\rangle  = \frac{1}{\epsilon^2}
\int_0^t \langle e_a,E^0_s[\Lambda_sN^{ji}_sX_s]\rangle
\langle C^{\theta_0}(Y_s)e_a,\mathrm{d}Y_s\rangle \\
& + 
\int_0^t
\langle e_a, Q^{\theta_0}(Y_s)E^0_s[\Lambda_sN^{ji}_sX_s]\rangle
 \mathrm{d}s
+ \delta_{aj}\int_0^t q^{\theta_0}_{ji}(Y_s)\langle
 e_i,E^0_s[\Lambda_sX_s]\rangle \mathrm{d}s,
\end{split}
\end{equation*}
which is equivalent to (\ref{new3}). Here, we have used that
\begin{equation*}
\begin{split}
 \mathrm{d}H_t &=  \langle e_a, X_{t-}\rangle \mathrm{d}N^{ji}_t
+ N^{ji}_{t-}\langle e_a,\mathrm{d}X_t\rangle + 
\langle e_a,\Delta X_t \rangle \Delta N^{ji}_t
\\ &= N^{ji}_{t-}\langle e_a,\mathrm{d}X_t\rangle + 
\delta_{aj}\langle e_i,X_{t-}\rangle \langle e_j,\mathrm{d} X_t \rangle.
\end{split}
\end{equation*}
The rest is similar and so omitted. \hfill////
\\

\noindent
{\bf   Conflict of Interest : } The author states that there is no conflict
of interest.

\end{document}